\documentclass[a4paper,12pt]{article}
\usepackage{amsmath}
\usepackage{amsfonts}
\usepackage{amsthm}
\usepackage{amssymb}
\usepackage{yfonts}
\usepackage{natbib}
\title{Harmonic ultrafilters}
\author{Rudi Hirschfeld\\
University of Antwerp\\
\texttt{rudihirschfeld@hetnet.nl}}
\newtheorem{theorem}{Theorem}
\newtheorem{corollary}{Corollary}
\newtheorem{remark}[theorem]{Remark}
\newtheorem{fact}[theorem]{Fact}

\newcommand{\betaN}{\beta\mathbb{N}}
\newcommand{\N}{\mathbb{N}}

\newcommand{\D}{\mathcal{D}}
\newcommand{\calH}{\mathcal{H}}

\numberwithin{equation}{section}
\begin{document}
\maketitle

\begin{abstract} A set of natural numbers will be called \emph{harmonic} if the reciprocals of its elements form a divergent series. An ultrafilter of the natural numbers will he called \emph{harmonic} if all each members are harmonic sets. The harmonic ultrafilters are shown to constitute a compact semigroup under the Glazer addition and its smallest ideal is obtained.

This paper is an extension of work begun by Hindman. Alle ingredients are found in the treatise by Hindman and Strauss.
\end{abstract}

\section {Introduction}
A set $A=\{a_1,a_2,\ldots\} \subset \N$ will be called \emph{harmonic} if the reciprocals of its terms form a divergent series,
\[\sum_{a_n\in A}\frac{1}{a_n}=\infty.\]
Sets with divergent reciprocals also are known as harmonically large and their counterparts as harmonically small.
Well-known examples of harmonic sets are $\N$ itself\footnote{The divergence of $\sum 1/n$ was already known to Nicole Oresme (c.1323 - 1382); he used an argument that is still taught in calculus 101 today.} and the set of the primes (known to Gauss; it is a direct consequence of his \emph{Primzahlsatz}). The main question concerning harmonic sets is the celebrated 

\vspace{4mm}

\noindent \textbf{Erd\"{o}s Conjecture}: \emph{does every harmonic set contain arithmetic progressions of arbitrary length?} 

This conjecture, which is inspired by the van der Waerden Theorem on arithmetic progressions in cells on integers, is and remains open. 

The family of all harmonic subsets of $\N$ will be denoted by $\calH$ and we start with a simple a Ramsey-type property of $\cal{H}$.
\begin{fact}\label{partition regular}
The set $\mathcal{H}$ is partition regular.
\end{fact}
\noindent \emph{Proof.} This means (see \cite{HinStr}, Definition 3.10 \footnote{We shall adhere to the terminology \emph{partition regular}, because it is short and standard, but it is a linguistic monstrum.}): if the union $\bigcup_{i=1}^r A_i$ of finitely many subsets $A_i$ of $\N$ belongs to $\mathcal{H}$, then there exists $i \in \{1,2,\ldots,r\}$ and $B\in \mathcal{H}$ such that $B\subset A_i$. As a matter of fact, there must be $i$ for which $A_i$ itself belongs to $\mathcal{H}$. If all $A_i$ were anharmonic, their union would be anharmonic too.\qed

Next,
\begin{fact}\label{translation-invariant}
The family $\mathcal{H}$ is left translation-invariant.
\end{fact}
\noindent \emph{Proof.} We have to show that: if $A\in \mathcal{H}$ and $s\in \N$, then $s+A\in \mathcal{H}$.
This invariance under left translations is immediate:
given $s\in \N$, put $N(s)=\min \{n\in \N: a_n\geq s\}$; then, for $N>N(s)$ we have
\[\sum_{n=1}^N\frac{1}{a_n+s}=\big(\sum_{n=1}^{N(s)}+\sum_{n=
N(s)+1}^N \big)\frac{1}{a_n+s}\geq
\sum_{n=1}^{N(s)}\frac{1}{a_n+s}+\sum_{n=N(s)+1}^N\frac{1}{2a_n}\]
and for $A\in \cal H$
the last sum tends to infinity together with $N$.\qed

\vspace{4mm}

We take some time out for a short commercial about Stone-\v{C}ech Compactifications. An extensive survey of the ideal theory in the compact semigroup $\betaN$ is found in \textsc{Hindman and Strauss} \cite{HinStr}. We only summarize some needed background information here.

As a set the Stone-\v{C}ech compactification $\betaN$ consists of all ultrafilters on $\N$. Each point of $\N$ itself is identified with the ultrafilter of its (discrete) neighborhoods. Next to these socalled fixed ultrafilters there must be a lot of \emph{free} ultrafilters $p$ (of one believes the Axiom of Choice in some form or another). For the subsets $A$ of $\N$ the \emph{closures} $\bar{A}=\{p \in \betaN: A \in p\}$ form a basis for the closed sets of a topology on $\betaN$ (but, strangely enough, also for the open sets). This way, $\betaN$ turns out to be a compact space and every map from $\N$ into a compact space has a unique continuous extension to all of $\betaN$ into the same compact space.
Because $(\N,+)$ is not a group, but a semigroup, shifting some $A \subset \N$ over $-x$, where $x\in\N$, must be defined with some care. One lets $A-x=\{y\in \N:y+x \in A\}$. Ordinary addition of integers is extended to $\betaN$ as follows: for two points $p$ and $q$ of $\betaN$, put
\[p+q\stackrel{def}{=}\{A \subset \N:\{x\in\N:A-x\in p\}\in q\}.\]
The extended addition is still associative and under this \emph{Glazer addition} $(\betaN,+)$ becomes a semigroup. Due to the asymmetry of $p$ and $q$ in this definition, this semigroup turns out to be highly non-commutative and addition is only continuous from the right.

\vspace{4mm}

A point $p\in \betaN$  will be called a \emph{harmonic ultrafilter} on $\N$ is all members $A$ of $p$ are harmonic sets.
The focus of this paper is directed to the set of all harmonic ultrafilters. 

In order to throw some light on the Erd\"{o}s Conjecture, \textsc{Hindman} in \cite{Hin} introduced the collection of all $p \in \betaN$ with the \emph{divergence} property that
\[A\in p \mbox{ if and only if } \sum_{a \in A} \frac{1}{a}=\infty.\]
In the present terminology $\D$ is the set of all harmonic ultrafilters on $\N$. It is shown in \cite{Hin} that 

\begin{fact}
$\D$ is a right ideal in $\betaN$. 
\end{fact}
For the sake of completeness we repeat Hindman's

\noindent \emph{Proof.} We have to show that $\mathcal{D}+\betaN \subset \mathcal{D}$. Take $p$ in $\mathcal{D}$ and $q$ somewhere in $\betaN$. Let $A\in p+q$. Will $\sum\{a^{-1}:a\in A\}=\infty$?.
Saying that $A$ belongs to $p+q$ means that we first have to consider those points $x\in \N$ for which $A-x \in p$, or $\sum\{b^{-1}:b\in A-x\}=\infty$. For the summand $b\in A-x$ there must exist $a\in A$ with $a=b+x$ (so that automatically $a>x$). In formula, $\sum\{(a-x)^{-1}: a\in A \mbox{ and }a>x\}=\infty$. This implies that
\begin{eqnarray*}
\sum\{\frac{1}{a}:a\in A\}&\geq& \sum\{\frac{1}{a}:a\in A \mbox{ and } a>x\}\\
&=&\sum\{\frac{1}{a-x}-\frac{x}{a(a-1)}:a\in A \mbox{ and }a>x\}\\&=&\sum\{\frac{1}{a-x}:a\in A \mbox{ and }a>x\}-x\sum\{\frac{1}{a(a-1)}:a\in A \mbox{ and }a>x\}.
\end{eqnarray*}
Whereas the second sum converges, the first one does not  
and we infer that every $A$ in $p+q$ is harmonic, indeed. \qed

\vspace{4mm}

In order to improve this result we first invoke some notation from Lemma 6.75 of \cite{HinStr}: \emph{the collection of those ultrafilters on $\N$ all of whose members are supersets of sets in $\mathcal{H}$ is designated by  $\Delta_{\mathcal{H}}$}. 

\begin{fact}
$\Delta_{\mathcal{H}} =\mathcal{D}$. 
\end{fact}

In fact, let $A \subsetneq \N$ belong to any $p\in \Delta_{\calH}$.  Take $b \in \N \setminus A$ and consider the proper superset $B=A\cup b$ of $A$. Since $A$ belongs to $\Delta_{\calH}$, all its supersets, in particular $B$, must be in $\calH$. Now, if $A\in p$ would not belong to $\calH$, then addition of the single point $b$ would turn it into the harmonic set $B$, which is impossible.\footnote{The equality $\Delta_{\calH}=\D$ means that if \emph{all} supersets of $A\subset \N$ are harmonic, then so is $A$ itself. For subsets it is the other way around, \emph{viz. every harmonic set contains an infinite anharmonic set}. Consider any harmonic set $A=\{a_1,a_2,\ldots\}$, together with an arbitrary anharmonic set $B=\{b_1,b_2,\ldots\}$. For each $k\in \N$ let $n_k$ be the smallest index with $a_{n_k}>b_k$. Write $c_k=a_{n_k}$. Then $C=\{c_1,c_2,\ldots\} \subset A$ and $\sum 1/c_k<\sum 1/b_k < \infty$.}

\begin{theorem}\label{D is 2-sided} The harmonic ultrafilters $p$ on $\N$ form a closed two-sided ideal $\mathcal{D}$ in $(\betaN,+)$.
\end{theorem}
\noindent\emph{Proof.} We shall invoke  Lemma 6.75 of \cite{HinStr}. The following three remarks are relevant here.
(1)
According to Fact \ref{partition regular} above, the family $\mathcal{H}$ of the harmonic subsets of $\N$ is partition regular. 
(2) The family  $\Delta_{\mathcal{H}}$ simply is $\mathcal{D}$. 
In fact, let $A \subsetneq \N$ belong to any $p\in \Delta_{\calH}$.  Take $b \in \complement A$ and consider the proper superset $B=A\cup b$ of $A$. Since $A$ belongs to $\Delta_{\calH}$, all its supersets, in particular $B$, must be in $\calH$. Now, if $A\in p$ would not belong to $\calH$, then addition of the single point $b$ would turn it into the harmonic set $B$, which is impossible.(3) The underlying collection $\calH$ is translation-invariant by the above Fact \ref{translation-invariant}. 
The Lemma 6.75 guarantees that $\mathcal{D}=\Delta_{\mathcal{H}}$ is a closed left ideal of $\betaN$. Together with the above result that $\D$ is a right ideal this implies that $\mathcal{D}$ is a two-sided ideal.\qed

The equality $\Delta_{\calH}=\D$ has a nice consequence:
\begin{corollary}
Every harmonic subset $A$ of $\betaN$ belongs to some point of $\D$; as a matter of fact, its closure $\bar{A}  \subset \D$. 
\end{corollary}
\noindent \emph{Proof.} Let $A\subset \N$ and $p \in \bar{A}$ (that is $A \in p$). For $A\in \calH$ all supersets $B$ trivially are harmonic. Hence $p\in \Delta_{\calH}$ by the definition of $\Delta_{\calH}$. Theorem \ref{D is 2-sided} guarantees that $p \in \D$.\qed

\vspace{4mm}

Being a closed ideal in $\betaN$, $\D$ is a compact semigroup for its own sake. By virtue of some fundamental results of Ellis, $\D$ has a unique minimal two-sided ideal $K(\D)$, generated \`a la $K(\D)=\D e\D$ by a minimal idempotent element $e \in \D$; \emph{cf.} Thm. 2.5, Thm. 2.9 and Lemma 1.49 in \cite{HinStr}. $K(\D)$ is called the \emph{smallest ideal} of $\D$; similarly, $\betaN$ has the non-empty smallest ideal $K(\betaN)$. 
\begin{theorem}
\[K(\D)=K(\betaN).\]
\end{theorem}
\noindent \emph{Proof.} This is a consequence of \cite{HinStr}, Thm.1.65 once we have shown that $K(\betaN) \cap \D \neq \emptyset$. But $\D$ is an ideal of $\betaN$ and as such it must contain the smallest ideal $K(\betaN)$; see \cite{HinStr}, Lemma 1.49. \qed

Sets $A \subset \N$ belonging to some $p\in K(\betaN)$ are called \emph{piecewise syndetic} sets. We recall that in terms of $\N$ itself, $A$ is piecewise syndetic if and only if the \emph{gaps} between its intervals of consecutive elements remain bounded in lengths, (see \cite{HinStr} Excercise 4.4.24). The above theorem yields
\begin{remark} If $A\in  \N$, then $A\in p$ with $p\in K(\D)$ if and only if $A$ is piecewise syndetic.
\end{remark}
A well-known version of the van der Waerden Theorem reads: every piecewise syndetic set contains arithmetic progressions of any length; see \cite{HinStr}, Thm. 14.1.  Such piecewise syndetic sets certainly are harmonic. On the other hand, the most famous harmonic set, to wit, the set of all primes, is not piecewise syndetic but is does contain arbitrarily long arithmetic progression, as has been shown by \textsc{Green and Tao} \cite{GreTao}.

\end{document}